\documentclass[12pt]{article}
\usepackage{amssymb,amscd,amsmath,amsthm}
\usepackage{latexsym,amstext}
\usepackage{latexsym,amstext}
\usepackage{color}
\usepackage{graphicx}
\usepackage{epstopdf}
\usepackage{cite}

\begin{document}

\title{\bf Approximate solutions of one dimensional systems with fractional derivative}

\author{A. Ferrari$^1$, M. Gadella$^2$, L.P. Lara$^3$, E. Santillan Marcus$^4$\\ \\
$^1$ Departamento de Matem\'atica, CONICET-Universidad Nacional de Rosario \\  Av. Pellegrini 250, 2000 Rosario, Argentina.\\[2ex]
$^2$ Departamento de F\'{\i}sica Te\'orica, At\'omica y Optica  and IMUVA, \\
Universidad de Va\-lladolid, 47011 Valladolid, Spain.\\[2ex]  $^3$Instituto de F\'isica Rosario, CONICET-UNR, \\ 
Bv. 27 de Febrero, S2000EKF Rosario, Santa Fe,  Argentina \\[2ex] $^4$Departamento de Matem\'atica, Universidad Nacional de Rosario, \\ Av. Pellegrini 250, 2000 Rosario, Argentina.
}

\maketitle

\begin{abstract}

The fractional calculus is useful to model non-local phenomena. We construct a method to evaluate the fractional Caputo derivative by means of a simple explicit quadratic segmentary interpolation. This method yields to numerical resolution of ordinary fractional differential equations. Due to the non-locality of the fractional  derivative, we may establish an equivalence between fractional oscillators and ordinary  oscillators with a dissipative term. 

\end{abstract}

\section{Introduction}

The study of fractional derivatives for its application in classical and quantum physics has lately received a lot of attention \cite{KST,HER,UCH}. Needless to say that one of the simplest and most studied of those systems is the one dimensional harmonic oscillator. Thus, it would be a good  point of departure in the study of systems with fractional derivative, a task which has been carried out in \cite{MAI}. Damped oscillator with fractional derivative has been also the objective of some studies, see \cite{OO}. Some extensions of the theory to other classical systems have been proposed, see for instance \cite{OO1} and references therein, or in \cite{CEY}. 

In many of these papers, it was noted an analogy between a fractional oscillator and a classical oscillator with a damping term. This could be an idea to be exploited in order to model quantum systems with dissipation, in which the second derivative of the wave function in the Schr\"odinger equation be replaced by a fractional derivative, see \cite{OO2}.

The present work has been inspired by the article by Narahari et. al. \cite{AHEC}, where in addition to the study of the one dimensional harmonic oscillator with fractional derivative, they give a comparison with an equivalent dissipative oscillator described on the phase plane and analyze the stability of the solution.

Along the present manuscript, we show that it may be possible the determination of a time interval in which the solution of a fractional one dimensional oscillator may be approximated by the solution of a one dimensional ordinary equation with a dissipative term. The idea could be described by using a very simple example. Let us consider the Caputo derivative $D^\alpha_0$, defined in \eqref{1} below, and the fractional differential equation $D^\alpha_0\,x(t)=0$, with $1<\alpha\le 2$ and initial conditions $x(0)=0$ and $\dot x(0)=-1$. The solution is $x(t)=-t$. Then, let us consider the equation $\ddot z(t)=-p\dot z(t)$, $p>0$. The goal is the determination of a value of $p$ such that the solution of this equation with initial conditions $z(0)=0$ and $\dot z(0)=-1$, i.e., the same initial conditions imposed to the fractional equation, be approximated by the solution of the fractional equation over a finite time interval. This is clear, since the solution of the equation on $z(t)$ is

$$
z(t)=\frac 1p (-1+e^{-pt})\,.
$$

Therefore, on a time interval $0<t<\tau$, with $\tau=1/p$, we have $z(t)=-t+o(t^2)$. In this sense of having similar approximate solutions on a finite interval, we say that the fractional and the dissipative equations are {\it equivalent}.  Here, we want to extend this idea. 

Observe that in our notion of equivalence, we have discarded the asymptotic regime. This is essentially due to two reasons: i.) for large values of time, the fractional oscillator does not show oscillations; ii.) the behaviour of the oscillator from a strictly physical point of view, whether linear or non-linear but particularly the latter, has interest for finite times only. Its asymptotic behaviour is not measurable and has a mathematical interest only, and it is not the object of our study.

The present paper is organized as follows: On Section 2, we construct a method to obtain approximate solutions of fractional differential equations, with fractional Caputo derivative to be defined there, based on segmentary interpolation. This kind of interpolation has been used successfully to obtain approximate solutions of ordinary differential equations \cite{GL}. On Section 3, we apply this method to the fractional linear oscillator and to some other simple examples and make estimations on its precision. We compare results with those obtained replacing the fractional oscillator by the ordinary oscillator with a dissipative term.  We present a similar analysis by replacing the equation of the oscillator by the van der Pol equation on Section 4. We close this paper with some concluding remarks.

\section{Caputo fractional derivative and its evaluation by segmentary interpolation}

Let $\alpha$ be a real positive number and denote by $n=\lceil \alpha \rceil$ the smaller integer bigger than $\alpha$. Let us define the Caputo fractional derivative, $D^\alpha_a$, of a $n$ times differentiable function of real variable, $x(t)$, as \cite{DIE}

\begin{equation}\label{1}
D^\alpha_a \,x(t) = \frac 1{\Gamma(n-\alpha)} \int_a^t \frac{x^{(n)}(s)}{(t-s)^{\alpha-n+1}} \,ds\,,
\end{equation}
where $x^{(n)}(s)$ means the $n$-th derivative of the function $x(s)$. Our objective, as mentioned at the header of the present section, is an evaluation of \eqref{1} using segmentary interpolation. Here, we consider that $0<\alpha<1$, so that the only choice for $n$ is $n=1$, and this will be the case for some of our applications. Segmentary interpolation is a standard tool of wide use in the approximation of solutions of differential equations \cite{HNW}. Let us sketch the method here for completeness, using an approach that has been used in previous articles by our group \cite{GL,FLS}.

Let $[a,b]$ be a compact interval in the real axis $\mathbb R$. At regular intervals, we select $n$ nodes, $a= t_0<t_1<\dots<t_n=b$, with $t_k-t_{k-1}=h$, for all $k=1,2,\dots,n$, so that $kn=b-a$. Let $x(t):[a,b] \longmapsto\mathbb C$ be a continuous function and use the notation $x_k:=x(t_k)$ and $I_k:=[t_{k-1},t_k]$\,, for all $k=1,2,\dots,n$. 

Then, a quadratic segmentary interpolator $S(t)$ for the function $x(t)$, is a continuous function $S(t):[a,b]\longmapsto \mathbb C$, with first continuous derivative, such that

1.- On each interval of the form $I_k=[t_{k-1},t_k]$, $k=1,2,\dots,n$, we have that $S(t)\equiv P_k(t)$, where $P_k(t)$ is a polynomial of order two, depending on the given interval. 

2.- The function $S(t)$ interpolates $x(t)$, in the sense that for any of the nodes $\{t_k\}$, one has that

\begin{equation}\label{2}
P_k(t_{k-1})=x_{k-1}\,,\qquad P_k(t_k)=x_k\,, \qquad k=1,2,\dots,n\,.
\end{equation}

The condition on the continuity of the derivative $S'(t)$ implies that

\begin{equation}\label{3}
P'_k(t_k) =P'_{k+1}(t_k)\,, \qquad k=1,2,\dots,n-1\,.
\end{equation}

Thus, the construction of the segmentary interpolator $S(t)$ relies in the construction of the interpolating polynomials $P_k(t)$. We propose the following form for the interpolating polynomials: For each of the intervals $I_k$, let us define,

\begin{equation}\label{4}
P_k(t) =p_k(t)+a_k(t-t_{k-1})(t-t_k)\,, 
\end{equation}
with

\begin{equation}\label{5}
p_k(t) =\frac{t-t_{k-1}}{h}\,x_k-\frac{t-t_k}{h}\,x_{k-1}\,,
\end{equation}
and the complex coefficients $a_k$ are given by

\begin{equation}\label{6}
a_k= \sum_{j=0}^n c_{k,j}\,x_j\,.
\end{equation}

We still need to determine the values of the $c_{k,j}$, which are

\begin{equation}\label{7}
c_{j,k} = \left\{ \begin{array}{ll} \frac{(-1)^k}{h^2}\,\eta_1\,,  & {\rm if}\; j=0\,, \\[2ex] 
\frac{(-1)^{k+1}}{h^2}\,(2\eta_1+\eta_2)\,, &  {\rm if}\; j=1\,, \\[2ex] 
\frac{(-1)^{k+j}}{h^2}\,(\eta_{j-1}+2\eta_j+\eta_{j+1}) \,, & {\rm if}\; 1<j<n-1\,, \\[2ex] 
\frac{(-1)^{k+n-1}}{h^2}\,(2\eta_{n-2}+\eta_{n-1}) \,, & {\rm if}\; j=n-1\,, \\[2ex] 
\frac{(-1)^{k+n}}{h^2}\, \eta_{n-1} \,, & {\rm if}\; j=n \,, 
\end{array}  \right.
\end{equation}
where $\eta_j=j/n$ if $j\le k-1$ and $\eta_j=j/n-1$ if $j>k-1$. 

Taking into account that $S(t)$ is an approximation of $x(t)$, on each of the nodes $t_k$ the Caputo fractional derivative \eqref{1} is approximated by

\begin{equation}\label{8}
D^\alpha_{a} x(t_k) \approx \frac{1}{\Gamma(1-\alpha)} \sum_{j=1}^k \int_{t_{j-1}}^{t_j} \frac{P'_j(s)}{(t_k-s)^\alpha}\,ds\,,\qquad k=1,\dots,n\,.
\end{equation}

Then, on each of the intervals $I_k$, we have that $P'_k(t)=\alpha_k\,t+\beta_k$, with

\begin{equation}\label{9}
\alpha_k=2a_k\,,\qquad \beta_k= \frac{x_k-x_{k-1}}{h} - a_k(2t_{k-1}+h)\,,
\end{equation}
and, consequently, equation \eqref{8}, takes the following form:

\begin{equation}\label{10}
D^\alpha_{a} x(t_k) \approx \frac{1}{\Gamma(1-\alpha)} \sum_{j=1}^k  \widetilde c_{k,j}\,\alpha_j + \widetilde d_{k,j}\,\beta_j\,.
\end{equation}

The new coefficients $\widetilde c_{k,j}$ and $\widetilde d_{k,j}$ are given by 

\begin{equation}\label{11}
\widetilde c_{k,j} = \int_{t_{j-1}}^{t_j} \frac{s\,ds}{(t_k-s)^\alpha}\,,\qquad  \widetilde d_{k,j} = \int_{t_{j-1}}^{t_j} \frac{ds}{(t_k-s)^\alpha}\,,
\end{equation}
which obviously depend solely of the partition. 

It is customary to choose $a=x_0=0$, which obviously does not restrict generality. Since the integrals in \eqref{11} are easily solvable and we know expressions for $\alpha_j$ and $\beta_j$, we can write the right hand side in \eqref{10} as

\begin{equation}\label{12}
D^\alpha_{a} x(t_k) \approx \frac{1}{(-1+\alpha)(-2+\alpha)\,\Gamma(1-\alpha)} \, \sum_{j=1}^k \gamma_{k,j}\,\alpha_j\,,
\end{equation}
where,

\begin{eqnarray}\label{13}
\gamma_{k,j} = [h(-j+k)]^{-\alpha} \,[h(1-j+k)]^{-\alpha} \, \Big\{ h(-1+j-k) \,[h(-j+k)]^\alpha \nonumber\\[2ex] \times
[-2h(-2+j+k) +h (-3+2j)\alpha-2(-2+\alpha)t_{j-1}] \nonumber\\[2ex] -h(j-k)\,[h(1-j+k)]^\alpha \, [-2h(-1+j+k)+h(-1+2j)\alpha-2(-2+\alpha) t_{j-1}] \Big\}\,.
\end{eqnarray}

This is a quite simple and workable receipt to obtain, once $x(t)$ is given, the values of its Caputo fractional derivative at the nodes $t_k$, so that we have an estimation of this fractional derivative. 

It is interesting to remark that, due to the linear dependence on $\{x_n\}$ of the coefficients $\alpha_j=2a_j$ given in \eqref{6}, then, the derivative $D^\alpha_{a} x(t_k)$ in \eqref{12} can be explicitly determined from $x(t)$.

\subsection{A type of differential equations with fractional derivative}

Let $x(t):[a,b]\longmapsto \mathbb R^m$ be a differentiable real function of the real variable $t$ and $f(t,x):[a,b]\times\mathbb R^m\longmapsto \mathbb R^m$. In addition, we assume that $x(t^*)=x^*$, where $t^*$ is one of the nodes $\{t_k\}$, $a\le t^*\le b$, and $0<\alpha<1$. Then, let us consider the following fractional differential equation:

\begin{equation}\label{14}
D^\alpha_a x(t) = f(t,x(t))\,.
\end{equation}

The objective is to obtain an approximation for the solution of equation \eqref{14} under the condition $x(t^*)=x^*$. We already know how to obtain the identity \eqref{14} in the nodes $t_k$. Take these nodes with the exception of $t^*$. Then, \eqref{14} provides of an algebraic system of equations where the indeterminates are $\{x_{j,k}\}_{j=1,k=0}^{m,n}$ with $x_{j,k}:=x_j(t_k)$ and $x_{j,k}\ne x^*_j=x_j(t^*)$. This algebraic system may or may not be linear depending on the form of $f(t,x(t))$ and is of order $(mn)\times(mn)$. A numerical solution of this system could be obtained by whatever method, which gives a segmentary solution $S(t)$, which is given once one has obtained the coefficients $a_k$ defined in \eqref{6}. 

In the particular case in which $f(t,x(t))$ contains an eigenvalue $\lambda$ and $f$ be linear with respect to $(\lambda,x)$ this algebraic system is linear and homogeneous. The eigenvalue is determined in the usual way. 

As the reader may easily understand, this method is more general than the usual way to obtain a solution knowing an initial value, since now $t^*$ could be any node. In particular, the restriction to the solution that replaces the initial value condition could be imposed at $t^*$, and this represents a great advantage when compared with the shooting method worked out in \cite{DIE1,DB}.

\section{The fractional linear oscillator}

A simple example of an equation of the type \eqref{14} is the linear oscillator with the fractional derivative, which is defined as

\begin{equation}\label{15}
D_0^\alpha\,x(t) =-\omega^2\,x(t)\,.
\end{equation}

As in the standard harmonic oscillator the constant $\omega^2=k/m$, where $m$ is the oscillator mass and $k$ a constant, $\alpha$ being the order of derivation that in the present case we assume to be $1<\alpha \le 2$. Using the definition \eqref{1}, taking into account that for some differentiable function $f(t)$ (in our case $f(t)=x(t)$ or $f(t)=\dot x(t)$, where the dot means first derivative), we have that

\begin{equation}\label{16}
\lim_{\alpha\to 0^+} \frac1{\Gamma(\alpha+1)}\, t^\alpha\,f(0)=f(0)\,,
\end{equation}
and that $n$ is either 2 or 3, we may integrate by parts \eqref{15} using \eqref{1}, which gives the following integral version of \eqref{2}:

\begin{equation}\label{17}
x(t) = x(0)+\dot x(0) -\frac{\omega^2}{\Gamma(-\alpha)}\int_0^t (t-s)^{-\alpha-1}\,x(s)\,ds\,.
\end{equation}

The general solution has the form

\begin{equation}\label{18}
x(t) = c_1\,  E_{\alpha,1}(-\omega^2 t^{\alpha}) + c_2\,t\, E_{\alpha,2} (-\omega^2 t^{\alpha})\,,
\end{equation}
where $E_{\alpha,\beta}(z)$ is the so called Mittag-Leffler function

\begin{equation}\label{19}
E_{\alpha,\beta}(z) = \sum_{k=0}^\infty \frac{z^k}{\Gamma(\alpha k+\beta)}\,.
\end{equation}

Thus, in order to obtain a particular solution, we have to impose some initial conditions. For instance, if we choose $x(0)=1$ and $\dot x(0)=0$, the solution to \eqref{17} with these initial conditions is given by

\begin{equation}\label{20}
x(t)= E_{\alpha,1} (-\omega^2 t^{\alpha})\,.
\end{equation}

Let us find a particular numerical solution to the fractional linear oscillator, using the method introduced in Section 2.1. We have to choose a particular value for $\omega$ and the simplest possibility is $\omega=1$. This is developed in the forthcoming subsection.

\subsection{Some numerical estimations.} 

First of all, it is not the objective here to give explicit expressions for the approximate solutions for the studied examples. It is not difficult to plot these solutions for different values of $n$.  

Let us start with equation \eqref{15} with $\omega=1$ on the interval $0\le t\le 1$, with $0<\alpha<1$ and initial condition $x(0)=1$. As seen above, this equation has exact solution given by $x_{\rm exact}(t)=E_{\alpha,1}(-t^\alpha)$ \cite{DIE}. The objective is now an estimation on the precision of the proposed method. As customary, this precision is measured by using the following parameter: 

\begin{equation}\label{21}
e_n(\alpha)= \int_0^1 [x_{\rm exact}(t)-x_n(t)]^2\,dt\,, 
\end{equation}

Here, $n$ is the number of sub-intervals $I_n$ in which we partite $[0,1]$, the number of nodes being $n+1$. The dependence of this parameter on $\alpha$ shows that the smaller is the value of $\alpha$, or equivalently the closer is $\alpha$ to zero, the lower is the precision and, therefore, the slower is the convergence to the exact value. However, we do not observe significative variations on the precision when we increase the value of $n$, i.e., as we make the sub-intervals smaller and smaller. 

\vskip1cm

\centerline{$
\begin{array}
[c]{cccc}
n & e_{n}(0.1) & e_{n}(0.5) & e_{n}(0.9)\\[2ex]
5 & 7.4\text{ }10^{-3} & 1.7\text{ }10^{-3} & 2.0\text{ }10^{-4}\\
10 & 3.0\text{ }10^{-3} & 4.9\text{ }10^{-4} & 2.8\text{ }10^{-5}\\
20 & 1.1\text{ }10^{-3} & 1.4\text{ }10^{-4} & 5.9\text{ }10^{-6}\\
40 & 3.0\text{ }10^{-4} & 3.7\text{ }10^{-5} & 1.3\text{ }10^{-6}
\end{array}
$}

\medskip
\centerline{Table 1.- Values of the precision in terms of $n$ and $\alpha$.}

\vskip1cm

This can be seen in Table 1, where we have chosen values of $n$ ranging from 5 to 40. The values of $\alpha$ studied are 0.1, 0.5 and 0.9. 

Let us study the precision of the method with another example different from the fractional oscillator, yet an equation of the form \eqref{14}. Here, we have chosen,

\begin{equation}\label{22}
D^{1/2}_0 x(t) =\sin x(t)\,,
\end{equation}
on the interval $[0,1]$, with the initial condition $x(1)=5/2$, which was already studied in \cite{DIE1}, where the integration was performed by means of the iterative shot method. Contrarily to the previous example, here we do not know an exact solution. The way out is to define the precision as

\begin{equation}\label{23}
e_n= \int_0^1 [D^{1/2}_0 x_n(t) -\sin x_n(t)]^2\, dx\,,
\end{equation}
where $n$ is again the number of intervals and $x_n(t)$ is the interpolating function for the studied case. After integration and using the boundary condition, we obtain $x(0)$. Along with \eqref{23}, we introduce another parameter that measures the convergence and that we denote as $e_r\%$. It represents the relative variation between the value of $y(0)$, obtained for a given value of $n$, and the value given for the precedent value of $n$ as listed on Table 2.  

 Table 2 is just a sample of numerous numerical examples we have performed. This sample is significative, as it manifest an obvious convergence and shows that the result obtained for a small number of nodes is satisfactory.  
 
\vskip1cm

\centerline{$
\begin{array}
[c]{cccc}
n & y(0) & e_{r}\% & e_{n}\\[2ex]
5 & 1.74895 & -- & 1.8\text{ }10^{-2}\\
10 & 1.73812 & 0.5 & 7.4\text{ }10^{-3}\\
20 & 1.73326 & 0.3 & 2.6\text{ }10^{-3}\\
30 & 1.73166 & 0.09 & 1.1\text{ }10^{-3}\\
40 & 1.73085 & 0.05 & 5.6\text{ }10^{-4}
\end{array}
$}
 
\medskip
\centerline{Table 2: Values of $y(0)$, $e_r\%$ and $e_n$ for a given value of $n$}

\vskip1cm

Finally, we have performed another reliability test, which was the use of the value of $y(0)$ obtained numerically as initial value and evaluate the value of $y(1)$. In all cases, we have recovered the value $y(1)=5/2$. 
 
\subsection{Damped oscillator with entire derivative.}

As is well known, the damped oscillator with entire derivative is given by

\begin{equation}\label{24}
m \ddot y(t)+ p\dot y(t)+k y(y)=0\,,
\end{equation}
where $m$, $p$ and $k$ are constants. Here, we assume that $p>0$. 

For $p=0$, the limit for $\alpha\longmapsto 2^-$ in \eqref{15} should give the solution $y(t)$ for \eqref{24}, which we denote here as $\lim_{\alpha\to 2^-}x(t)=y(t)$. The solution $x(t)$ is damped oscillatory on the transitory regime only \cite{AHEC,DIE,KST}. Based on these notions, we propose the following Ansatz: 

{\it For each given $1<\alpha\le 2$, there exists $p>0$  such that the solution $x(t)$ of \eqref{15} with $\alpha$ is a good approximation of the solution $y(t)$ of \eqref{24} with $p$, in the transitory regime}.

Using this Ansatz, let us obtain an approximate solution $y(t)$ for \eqref{24} such that this and its corresponding solution $x(t)$ for \eqref{15} fulfil the conditions $y(0)=x(0)$ and $\dot y(0)=\dot x(0)$. This is:

\begin{equation}\label{25}
y(t) = \exp\left( -\frac pm\,t \right)\, (c_-\exp(-\Delta t) +c_+ \exp(\Delta t))\,,
\end{equation}
where, 

\begin{equation}\label{26}
\Delta= \sqrt{(p/m)^2-4\omega^2}\,,\quad \lambda_\pm =\frac{-p\pm\Delta}{2}\,,\quad c_\pm =\pm \frac{\lambda_\mp}{\Delta}\,,
\end{equation}
where $\omega^2$ was given in \eqref{15}. 

Then, the point is the determination of the value of $p$ being given the value of $\alpha$, or equivalently the determination of a function $p=p(\alpha)$, in application to our Ansatz. This is an optimal control problem. We have to find the optimal solution, which minimizes  the following functional:

\begin{equation}\label{27}
E(\alpha):= \frac 1T \int_0^T [x(t)-y(t)]^2\,dt\,,
\end{equation}
where $T$ is some time scale in which the amplitude of the oscillations are reduced by a factor of $1/T$. On the interval $[0,T]$, the transitory regime, we compare the solutions of the fractional derivative $x(t)$ and of the damped oscillator $y(t)$.  The functional $E(\alpha)$ measures the deviation between $x(t)$ and $y(t)$. Then, go back to \eqref{20} and note that the function $E_{\alpha,1} (-\omega^2 t^{\alpha})$ is not asymptotically oscillating as $t\longmapsto 0$. This permits us to choose a value of $T$, although not small, not very high either. Numerical experiments have shown that the choice $T=20$ is appropriate. 

Let us give some numerical results. On Table 3, we give the dependence between values of $\alpha$, $p(\alpha)$ and $E(\alpha)$ for the values $k=m=1$ and $n=50$. 

\vskip1cm

\centerline{$
\begin{array}
[c]{ccc}
\alpha & p & E\\[2ex]
1.10 & 1.140 & 6.4\,10^{-3}\\
1.30 & 0.891 & 5.7\,10^{-3}\\
1.50 & 0.668 & 4.7\,10^{-3}\\
1.70 & 0.433 & 3.3\,10^{-3}\\
1.90 & 0.152 & 1.1\, 10{-3}\\
1.95 & 0.754 & 4.3\, 10^{-4}\\
2 & 0 & 0
\end{array}
$}

\medskip
\centerline{Table 3: Comparison between the values of $\alpha$, $p$ and $E$,}

\centerline{for $T=20$, $k=m=1$ and $n=50$.}

\vskip1cm 

An explicit expression of the function $p(\alpha)$ may be obtained by the least-square method and this gives $p(\alpha)= 1.49409 +0.056127 \,\alpha - 0.401446\,\alpha^2$. This is depicted on Figure 2.

On Figure 1, we represent the usual behaviour in the transitory regime for $x(t)$ and $y(t)$, when we choose $\alpha=1.7$ and $n=50$. 

It is interesting to remark that numerical experiments show that $p(\alpha)$ does not depend on any choice of initial values. 

\subsection{Non-linear oscillator}

Following the discussion on the damped oscillator, we present a similar problem given by the following non-linear oscillator:

\begin{equation}\label{28}
D_0^\alpha\,x(t) =y(t)\,, \qquad D_0^\alpha\,y(t) =-\sin x(t)\,,
\end{equation}
with $0<\alpha\le 1$. Let us choose the initial values given by $x(0)=1$ and $y(0)=0$. Clearly, for small oscillations equation \eqref{28} becomes equation \eqref{15} with the replacement $\alpha \rightarrow 2\alpha$. Again, this is a non-linear problem having no analytic solution for $\alpha$ non-integer. Then, we proceed by analogy with the damped oscillator. In this case, we consider the following system involving entire derivatives only:

\begin{equation}\label{29}
\dot z(t)=w(t)\,,\qquad \dot w(t)=-p w(t)-\sin z(t)\,,\qquad p>0\,.
\end{equation}

On the transitory regime, we compare the solutions of systems \eqref{28} and \eqref{29} under the conditions $z(0)=x(0)=1$ and $w(0)=y(0)=0$. To do this, we need the previous determination of $p(\alpha)$, which we assume that minimizes the following quadratic dispersion: 

\begin{equation}\label{30}
E(\alpha)=\frac 1T \int_0^T  \{ [x_n(t)-z(t)]^2 +[y_n(t)-w(t)]^2 \}\,dt\,.
\end{equation}

Obviously, this expression generalizes \eqref{27}. Again, we adjust the value of $T$ by numerical experiments, which show that $T=20$ is, again, a convenient choice. On Table 4, we give some values for the dependence between $\alpha$, $p(\alpha)$ and $E(\alpha)$ after the choice $T=20$ and $n=50$. 

\vskip1cm

\centerline{$
\begin{array}
[c]{ccc}
\alpha & p & E\\[2ex]
0.50 & 1.203 & 5.5\,10^{-3}\\
0.70 & 0.757 & 3.9\,10^{-3}\\
0.90 & 0.294 & 2.1\,10^{-3}\\
0.95 & 0.148 & 1.4\,10^{-3}\\
1.00 & 0 & 0
\end{array}
$}

\medskip
\centerline{Table 4: Comparison between the values of $\alpha$, $p$ and $E$,}

\centerline{for the values $T=20$ and $n=50$.}

\vskip1cm

The above examples manifest an analogous behaviour between a fractional  linear oscillator and a damped or even non-linear oscillator on some time interval. The solutions between the fractional and the integer equation are very similar on some time scale. This could be a rather general situation, so that in many practical cases and inside a time interval, we conjecture that a fractional operator might be replaced by the frictional additional term on the classical oscillator. The behaviour of the solutions is similar to that shown in Figure 1. 

\begin{figure}
\centering
\includegraphics[width=0.4\textwidth]{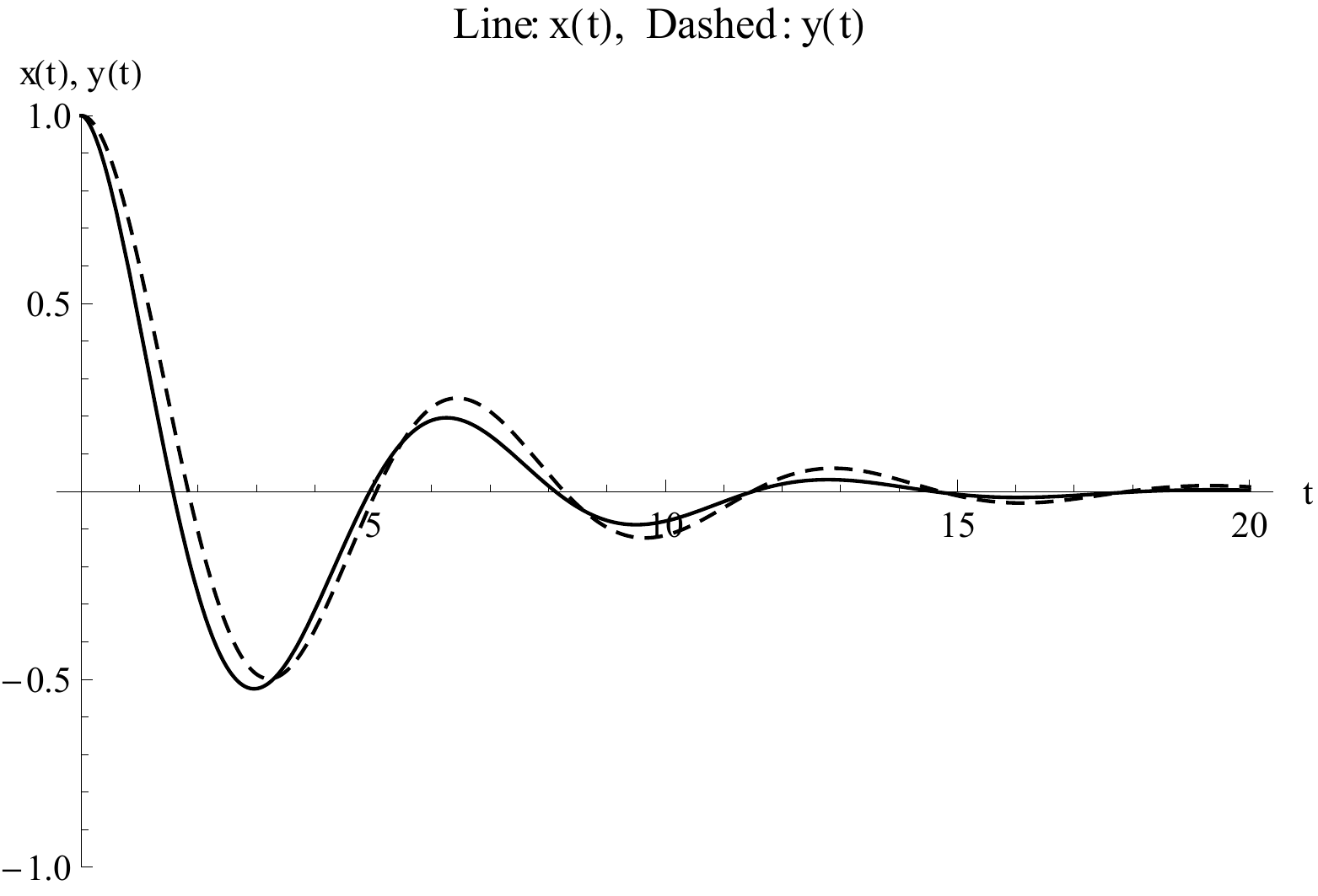}
\caption{\small  
{ The continuous line represents the solution, $x(t)$ of the fractional equation \eqref{15}, which is given by \eqref{25}. The dashed line gives the solution, $y(t)$, of the equation with ordinary derivative \eqref{25}.} 
}
\label{}
\end{figure}

\begin{figure}
\centering
\includegraphics[width=0.4\textwidth]{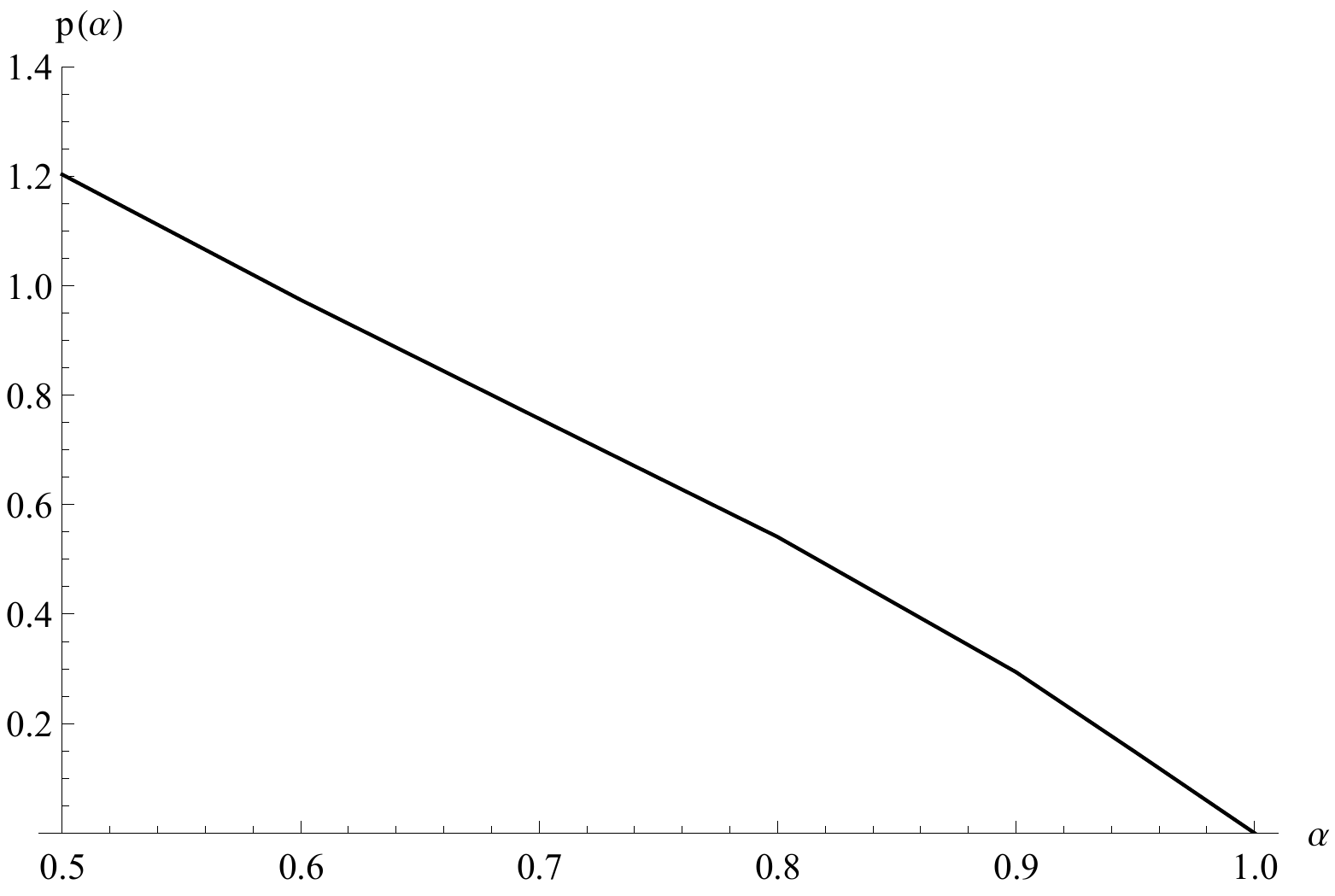}
\caption{\small  
Function $p(\alpha)$. 
}
\label{}
\end{figure}

\section{On the fractional van der Pol equation}

The van der Pol equation is a second order non-linear ordinary differential equation \cite{VDP,FAR}. It has the following form:

\begin{equation}\label{31}
\ddot x(t) - \mu(1-x^2(t))\,\dot x(t) + x(t)=0\,,
\end{equation}
where $\mu\ge 0$ is a constant. When $\mu=0$, \eqref{31} is the equation of the ordinary harmonic oscillator. The van der Pol equation may be written in terms of a first order system as

\begin{equation}\label{32}
 \dot x(t):= z(t)\,,\qquad  \dot z(t)=-\mu (x^2(t)-1) z(t) -x(t)\,, \qquad \mu\ge 0\,.
\end{equation}

This equation has a unique limit cycle for $\mu\ne 0$, after the Li\'enard theorem \cite{STR}. This suggested us the interest of considering the possible existence of a limit cycle for the fractional system analogous to \eqref{32} given by

\begin{equation}\label{33}
D_0^\alpha\, x(t)=z(t)\,,\qquad D_0^\alpha \,z(t)= -\mu z(t) (x^2(t)-1)-x(t)\,,\qquad 0<\alpha\le 1. 
\end{equation}

One fractional van der Pol equation of the type

\begin{equation}\label{34}
D_0^{\alpha+1}\, x(t) +\mu (x^2(t)-1)  D_0^\alpha \,x(t) +x(t)=0\,,
\end{equation}
has been studied in \cite{GUO}, where a relation between the parameters $\alpha$ and $\mu$ was given as a sufficient condition for the existence of a limit cycle, using the balance harmonic method \cite{DEU}.  

We have studied the system \eqref{33} through numerical as well as analytic methods. We have performed a big amount of numerical experiments, which have shown the existence of a value of the parameter $\mu$, here called $\mu_c$, where the subindex $c$ stands for {\it critical}, which depends on $\alpha$ and $\mu_c(\alpha)>0$, such that

\begin{itemize}

\item{For values of $\mu$ with $0<\mu<\mu_c$, there is a fixed point $(x^*,z^*)=(0,0)$, which remains stable at the limit $t\longmapsto\infty$, $\lim_{t\to\infty}(x(t),z(t))=(0,0)$. Therefore, there is no stable limit cycle. In addition, there is no evidence of the existence of an unstable limit cycle.}

\item{For values $\mu>\mu_c$, the fixed point $(x^*,z^*)=(0,0)$ is unstable. We found a unique stable limit cycle. In this case a Hopf bifurcation emerges with $\mu_c$ as critical parameter.}

\item{As shown in Figure 3, $\mu_c(\alpha)$ decreases with $\alpha$ and $\lim_{\alpha \to1}\mu_c(\alpha)=0$.}

\end{itemize}

The point here is to show  the existence of the critical value for the parameter $\mu$ for a given value of $0<\alpha<1$, $\mu_c(\alpha)$. This existence has been manifested by the numerical estimations above. Nevertheless, this existence may be also shown analytically. To this end, we use the following result \cite{AEE}:

Let us consider the following system, where $D_0^\alpha$ represents the Caputo fractional derivative: 

\begin{equation}\label{35}
D^\alpha_0 \, x(t)=f(x,z)\,,\qquad D^\alpha_0\,z(t)= g(x,z)\,,\qquad 0<\alpha<1\,.
\end{equation}

A solution $(x^*(t),z^*(t))$ is in equilibrium if $f(x^*(t),z^*(t))=g(x^*(t),z^*(t))=0$. It is asymptotically stable if the eigenvalues, $\lambda$, of the Jacobian matrix

\begin{equation}\label{36}
J(x,z):= \left(\begin{array}{cc} \partial f/\partial x & \partial f/\partial z \\[2ex] \partial g/\partial x & \partial g/\partial z  \end{array} \right)\,,
\end{equation}
when evaluated at the equilibrium point satisfies

\begin{equation}\label{37}
|\arg(\lambda)|>\alpha\,\frac \pi 2\,.
\end{equation}

A comparison between \eqref{33} and \eqref{35} gives the precise form of $f(x,z)$ and $g(x,z)$ for our particular case. This gives the precise form of \eqref{36} as

\begin{equation}\label{38}
J(x,z):= \left(\begin{array}{cc} 0 & 1 \\[2ex] -1-2\mu\,z(t)x(t) & -\mu(x^2(t)-1)  \end{array} \right)\,.
\end{equation}

Taking into account that the fixed point is located at $(x^*,z^*)=(0,0)$, we have that

\begin{equation}\label{39}
J(0,0)=\left(\begin{array}{cc} 0 & 1\\[2ex] -1 & \mu  \end{array} \right)\,,
\end{equation}
which has the following eigenvalues:

\begin{equation}\label{40}
\lambda_\pm= \frac 12 (\mu\pm \sqrt{\mu^2-4})\,.
\end{equation}

Obviously, if $\mu \ge 2$, then $\arg(\lambda_\pm)=0$. On the other hand, if $\mu<2$, one has that,

\begin{equation}\label{41}
\arg(\lambda_\pm)= \arctan \left(\pm\sqrt{\left( \frac 2\mu\right)^2-1} \right)\,.
\end{equation}

From \eqref{37}, the critical value, $\mu_c$, of $\mu$ should obey the following relation:

\begin{equation}\label{42}
\arctan \left(\pm\sqrt{\left( \frac 2\mu_c\right)^2-1} \right) = \alpha\,\frac \pi 2\,,
\end{equation}
which gives

\begin{equation}\label{43}
\mu_c= \frac{2}{\displaystyle \sqrt{1+\tan^2\left( \frac{\alpha\pi}{2}\right)}}\,.
\end{equation}

\begin{figure}
\centering
\includegraphics[width=0.4\textwidth]{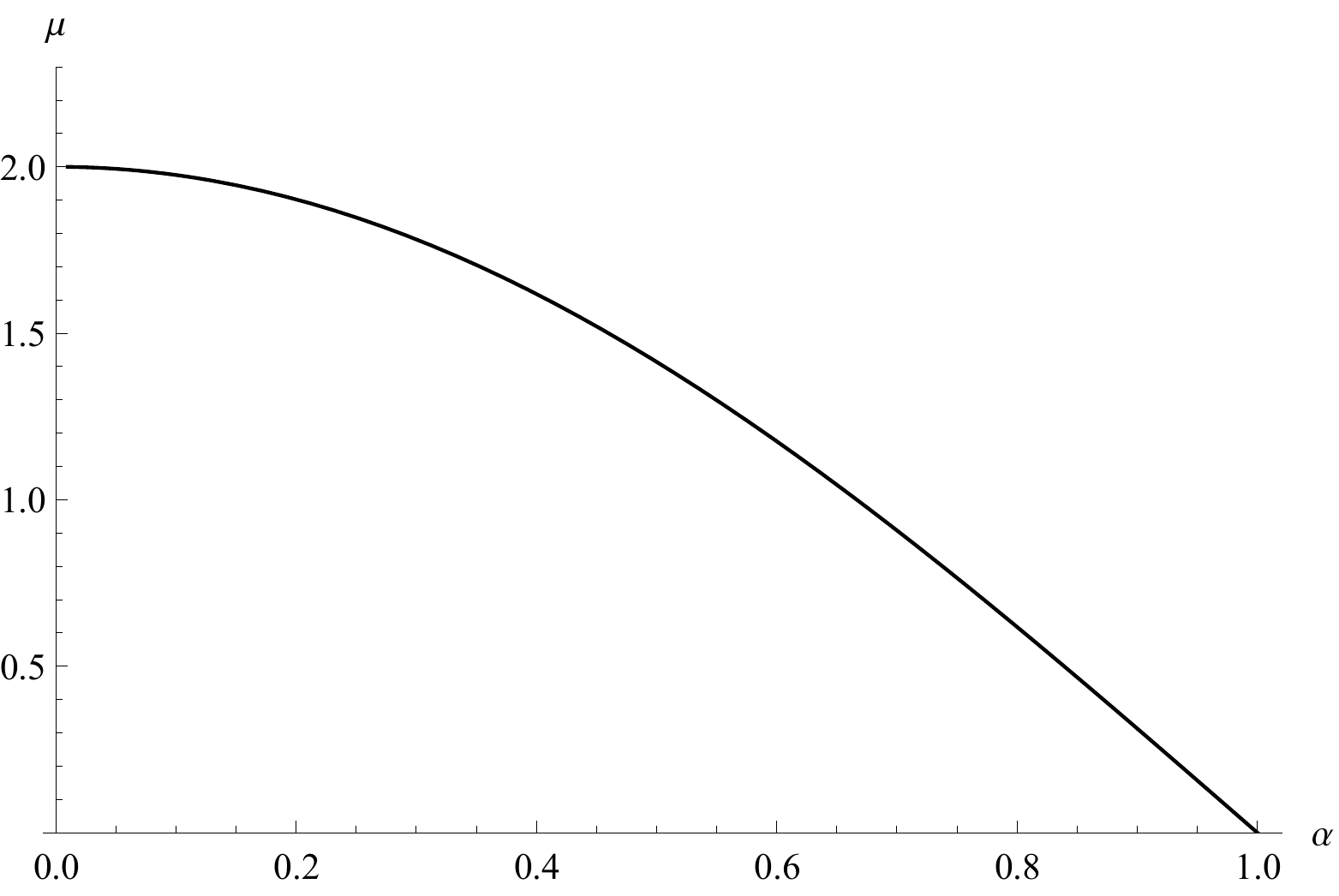}
\caption{\small  
Function $\mu_c(\alpha)$. 
}
\label{}
\end{figure}

This is to say, if we fix $\alpha$ and start from $\mu\approx 0$, as we increase $\mu$, we go from a situation with a asymptotically stable fixed point to an unstable point. This happens when $\mu>\mu_c$. The transition from the stability to the unstability drives to the emergency of a limit cycle. We may qualitatively interpret the limit cycle loss as follows: let us consider $\mu\approx 0$ in \eqref{33}, which may then be approximated  by

\begin{equation}\label{44}
D_0^\alpha \,x(t) = z(t)\,,\qquad D_0^\alpha \,z(t)=-x(t)\,.
\end{equation}

This is the same than equation \eqref{28} with the paraxial approximation $\sin y(x)\approx y(x)$. Note that \eqref{44} does not show a limit cycle and, further, the trivial solution $(0,0)$ is an attractor. The second equation in \eqref{33} contains the term $-\mu\,z(x)(y^2(x)-1)$, which in the case of $\mu>\mu_c$ is not negligible. This fact outbalances the dissipation and this is precisely which makes it possible the existence of a limit cycle.

On Figure 4 and on the phase plane, the continuous and slashed curves represent the solution with entire and fractional derivative, respectively. Both trajectories are determined with same initial values and same parameter $\mu$. In all cases, the fractional limit circle is enclosed by the trajectory of the limit cycle with entire derivative. 

On Figure 3, we show the relation $\mu_c=\mu_c(\alpha)$.  In the region above the curve, there exists a stable cycle limit and, furthermore, the fixed point $(0,0)$ is unstable. Below the curve the limit cycle does not exist and the fixed point is asymptotically stable. There is an obvious difference with the results obtained in \cite{GUO}, which is due to the fact that the fractional equations \eqref{33} and \eqref{34} are {\it not} equivalent.

\subsection{Equivalence between the fractional van der Pol equation and the same equation with entire derivative and dissipation.}

On the previous section, we have compared the approximate solutions of a dissipative oscillator with entire derivative with those of the linear oscillator with fractional derivative. Now, we want to carry out a similar analysis with the fractional van der Pol equation and a van der Pol equation with entire derivative and a dissipative term of the form $\beta z(t)$, $\beta>0$.   This system has the form,

\begin{equation}\label{45}
x'(t)=z(t)\,,\qquad z'(t)=-z(t)(\beta+\mu(x^2(t)-1))-z(t)\,, \qquad \mu\,,\beta>0\,.
\end{equation}

Here, the fixed point is $(x^*,z^*)=(0,0)$. To check its stability, we consider again equation \eqref{37}, which in the present case gives at the fixed point the following expression

\begin{equation}\label{46}
J(0,0)= \left(\begin{array}{cc} 0 & 1 \\[2ex] -1 & \mu-\beta \end{array}\right)\,,
\end{equation}
which has the eigenvalues

\begin{equation}\label{47}
\lambda_\pm = \frac{(\mu-\beta) \pm \sqrt{(\mu-\beta)^2-4}}2\,.
\end{equation}

Therefore, the fixed point is stable if Re$(\lambda_\pm)<0$ and unstable if Re$(\lambda_\pm)>0$, or equivalently, if $\mu-\beta<0$ and $\mu-\beta>0$, respectively. Then, for each $\mu$, there exists a $\beta_c=\mu$, where the Hopf bifurcation of the fixed point appears and, consequently, the destruction of the limit cycle. 

In any case, according to the Li\'enard theorem \cite{STR,FAR}, we may show that there exists a unique stable limit cycle if $\mu-\beta<0$. In consequence, the van der Pol equations with entire derivative and dissipation and fractional have qualitatively the same properties.

We have checked numerically a qualitative equivalence, in the sense of having approximately the same solution, through a substantial number of numerical experiments, between equations \eqref{44} (with fractional derivative) and \eqref{45} (with entire derivative). Thus by trial and error, we have determined a value of $\beta$ giving the same cycle in both cases. For instance, if we give the values $\alpha=0.9$ and $\mu=0.1$, we obtained $\beta\cong 0.315$. In an analogous manner, we trial with values for which $\mu<\mu_c$ and obtained similar conclusions.

On Figure 5, we represent limit cycles for the fractional and dissipative entire van der Pol equations.

\begin{figure}
\centering
\includegraphics[width=0.4\textwidth]{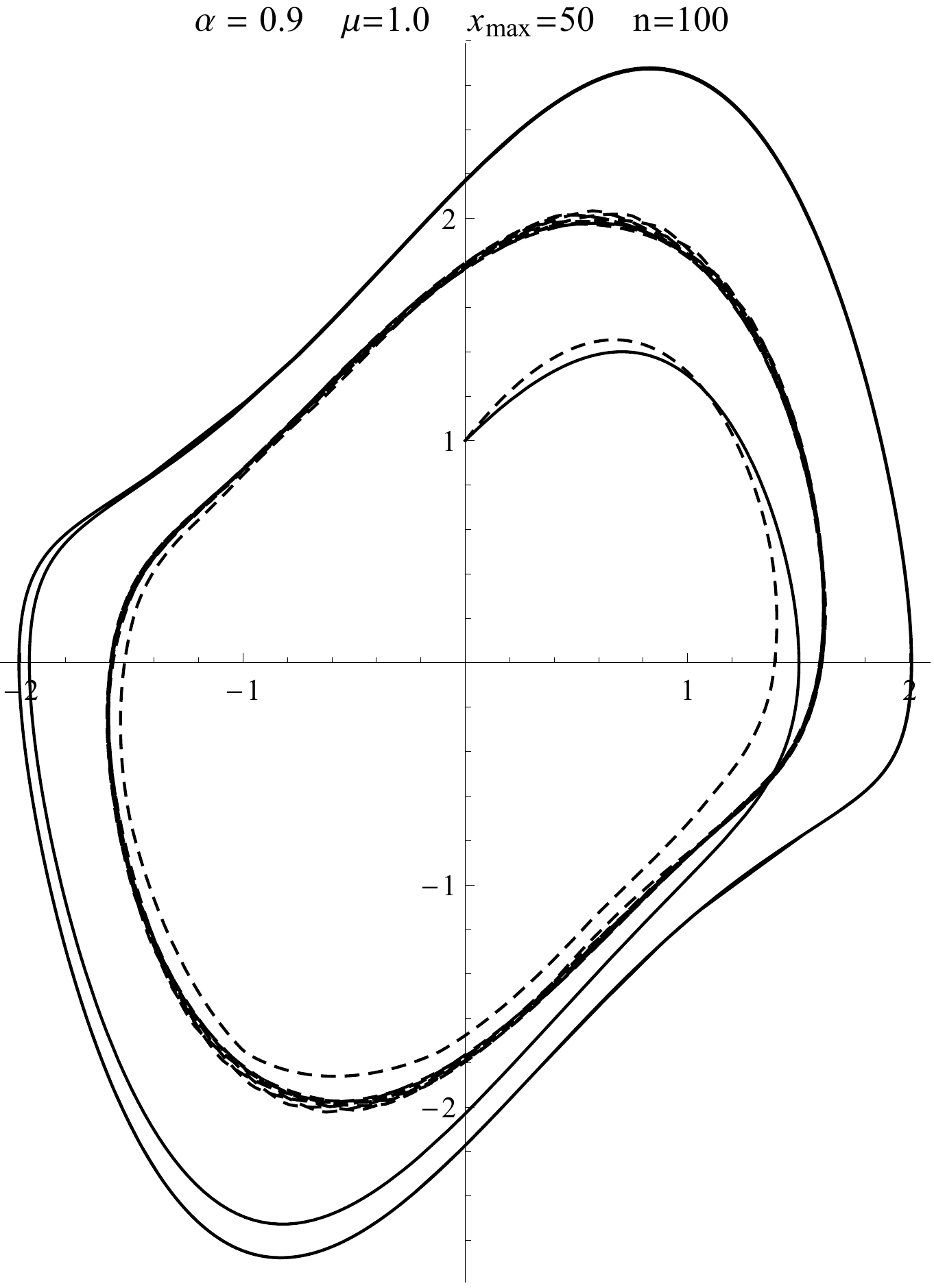}
\caption{\small  
Comparison between the entire (continuous curve) and the fractional (slashed curve) van der Pol solutions, with the same value of $\mu$. 
}
\label{}
\end{figure}

\begin{figure}
\centering
\includegraphics[width=0.4\textwidth]{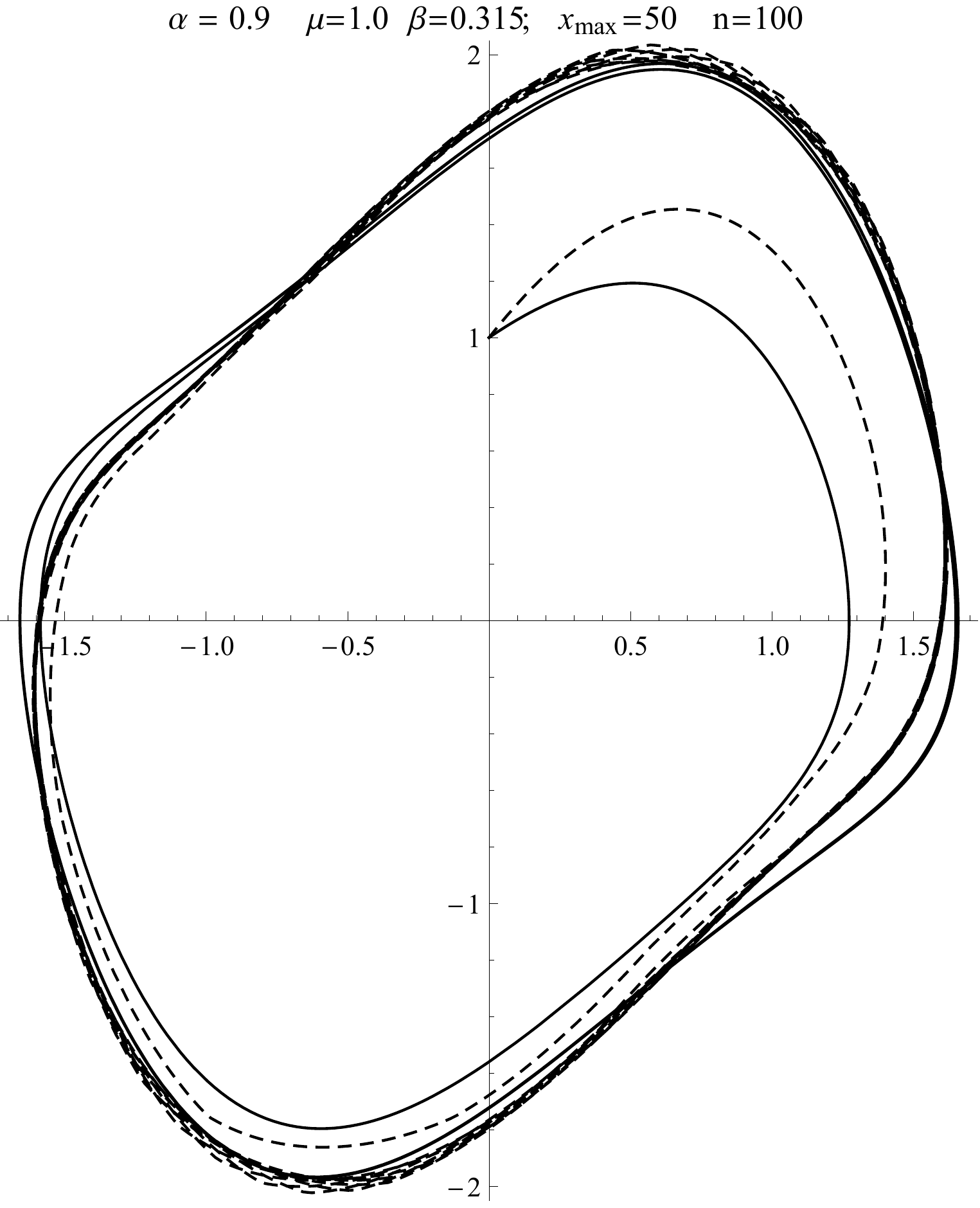}
\caption{\small  
Comparison between the damped (continuous curve) and the fractional (slashed curve) van der Pol solutions. 
}
\label{}
\end{figure}

\section{Concluding remarks}

We have applied a quadratic spline method  in order to obtain functions that approximate the result of applying a fractional derivative to a given function.  This is suitable to obtain solutions to some differential equations with initial values or mixed conditions of potential interest in Physics or Engineering.  We have tested our method with the fractional linear oscillator, where exact solutions are known and checked its degree of precision. Results of numerous numerical experiments show that for values of $\alpha$ in the range $0<\alpha<1$, the higher is $\alpha$, the better is the precision of our method.  Here, $\alpha$ is the order of the fractional derivation, $D^\alpha$. However, there are not substantial differences when we increase the number of nodes on the interval under our consideration. Similar results have been obtained for non-linear oscillators. 

Based on previous studies on the approximation of solutions of the fractional linear oscillator by solutions of a damped oscillator, we have used our method to confirm these results. We have shown that there exists a time interval for which the solutions of both equations are similar with a high degree of accuracy, under the condition that a relation is given between the coefficient $p$ of the dissipative term of the damped oscillator and the order of the fractional derivation, $\alpha$. 

A similar study compares a fractional and a damped van der Pol equations, written as a system of two equations on phase space, with similar results.  In addition, we have considered the behaviour of limit cycles and fixed points in terms of $\alpha$ and a parameter $\mu$ characteristic of the van der Pol equation. Using analytic as well as numerical arguments, we show the existence of a critical value for the parameter $\mu$, $\mu_c$, such that if $\mu<\mu_c$ the origin of phase space is stable and if $\mu>\mu_c$ is unestable.  This limit value $\mu_c$ depends on $\alpha$ and we give the exact relation.

\section*{Acknowledgements}  

This research has been financed by  the Projects No. ING 19/i 402 and ING 80020180100064  of the Universidad Nacional de Rosario, the Spanish MINECO (Project  No. MTM2014-57129) and the Junta de Castilla y Le\'on (Project Nos. BU229P18 and VA137G18).


\begin{thebibliography}{99}

\bibitem{KST} A. Kilbas, H. Srivastava, J. Trujillo. Theory and Applications of Fractional differential equations, Elsevier, Amsterdam (2006). 

\bibitem{HER} R. Herrmann,  {\it Fractional Calculus: An Introduction for Physicists}. World Scientific,  (2011).

\bibitem{UCH} V.V. Uchaikin, {\it Fractional derivatives for physicists and engineers.
Volume I Background and theory volume, Volume II Applications}. 
Berlin
Heidelberg: Springer Science \& Business Media (2013).

\bibitem{MAI} F. Mainardi, Fractional relaxation-oscillation and fractional diffusion-wave phenomena, Chaos, Solitons and Fractals, {\bf 7} (9), 1461-1477 (1996). 





\bibitem{OO} F. Olivar-Romero, O. Rosas-Ortiz, Fractional driven damped oscillator, J. Phys. Conf. Ser., {\bf 839}, 012010 (2017). 

\bibitem{OO1}  F. Olivar-Romero, O. Rosas-Ortiz, Transition from the Wave Equation to Either the Heat
or the Transport Equations through Fractional Differential Expressions, Symmetry, {\bf 10}, 524 (2018). 

\bibitem{OO2}  F. Olivar-Romero, O. Rosas-Ortiz, Factorization of the quantum fractional oscillator, J. Phys. Conf. Ser., {\bf 698}, 012025 (2016). 

\bibitem{CEY}  Can Evren Yarman,  Approximating fractional derivative of Faddeeva function, Gaussian function, and Dawson's integral, Math. Met. Appl. Scie., DOI: 10.1002/mma.5679 (2019). 





\bibitem{AHEC} B.N. Narahari Archar, J.W. Hanneken, T. Enck, T. Clarke, Dynamics of the fractional oscillator, Physica A, {\bf 297}, 361-367 (2001). 

\bibitem{DIE} K. Diethelm, {\it The Analysis of Fractional Differential Equations. An Application Oriented Exposition Using Differential Operators of Caputo Type}, Springer Verlag, Berlin (2010). 

\bibitem{DIE1} K. Diethelm, W. McLean, Volterra integral equations and fractional calculus: do neighboring solutions intersect?, Journal of Integral Equations and Applications, {\bf 24} (1), 25-37 (2012).

\bibitem{HNW} E. Hairer, S.P. N{\o}rsett and G. Wanner, {\it Solving Ordinary Differential Equations I}, Springer, Berlin and New York, 1993.

\bibitem{GL} M. Gadella, L.P. Lara, A study of periodic potentials based in quadratic splines, Int. J. Mod. Phys. C, {\bf 29}, 1850067 (2018).

\bibitem{FLS} A. Ferrari, L.P. Lara, E Santillan-Marcus, Convergence analysis and parity conservation of a new form of quadratic splines, arXiv:1906.10559v1 (2019). 







\bibitem{DB} H. Demir, Y. Balturk, On numerical solutions of fractional order boundary value problem with shooting method, ITM Web of Conferences, {\bf 13}, 01032 (2017). 

\bibitem{VDP} B. van der Pol, On relaxation equations, The London, Edinburgh and Dublin Phil. Mag. \& J. of Sci., 2 (7), 978-992 (1927).



\bibitem{STR} S. Strogatz, {\it Nonlinear Dynamics and Chaos}, CRC Press, Taylor and Francis, Boca Raton, London and New York, 2015. 

\bibitem{FAR} M. Farkas, {\it Periodic Motions}, Springer, New York and Berlin (1994).

\bibitem{GUO} Z. Guo, A.Y.T. Leung, H.X. Yang, Osillatory region and asymptotic solution of fractional van der Pol oscillator via residue harmonic balance technique, Apply Mathematical Modelling, {\bf 35}, 3918-3925 (2011). 

\bibitem{DEU} P. Deuflhard, {\it Newton Methods for Nonlinear Problems}, Springer, Berlin, 2006.



\bibitem{AEE} E. Ahmed, A. El-Sayed, A.A. El-Saka, On some Routh-Hurwitz conditions for fractional order differential equations and their applications in Lorenz, R\"ossler, Chua and Chen systems, Phys. Lett. A, {\bf 358}, 1-4 (2006). 


\end{thebibliography}
\end{document}